\documentclass{svmult-ddm}

\usepackage{mathptmx}       
\usepackage{helvet}         
\usepackage{courier}        
\usepackage{type1cm}        
\usepackage{graphicx}        

\usepackage[bottom]{footmisc}

\usepackage{amssymb}
\usepackage{amsfonts}
\usepackage{amsbsy}
\usepackage{amscd}
\usepackage{amstext}
\usepackage{dsfont}
\usepackage[english]{babel}
\usepackage{color}
\usepackage{graphics}
\usepackage{epsfig}
\usepackage{subfigure}
\usepackage{wrapfig}
\usepackage{psfrag}
\usepackage{color}
\usepackage{url}
\usepackage{verbatim}
\usepackage{algorithm}
\usepackage{algorithmic}

\begin{document}
\title*{Aggregation-based aggressive coarsening with polynomial smoothing}
 
  

\author{James~Brannick\inst{1}}
\institute{Department of Mathematics, The Pennsylvania State University,
University Park, PA 16802, USA
\texttt{brannick@psu.edu}}
%
%
\maketitle

\abstract{
This paper develops an algebraic multigrid preconditioner for the graph Laplacian.
The proposed approach uses aggressive coarsening based on the aggregation framework in the setup phase and a polynomial smoother with sufficiently large degree within a (nonlinear) Algebraic Multilevel Iteration as a preconditioner to 
the flexible Conjugate Gradient iteration in the solve phase. 
We show that by combining these techniques it is possible to
design a simple and scalable algorithm.  
Results of the algorithm applied to graph Laplacian systems arising 
from the standard linear finite element discretization of the scalar Poisson problem
are reported.
}

\section{Introduction}
This paper concerns the development of an algebraic multigrid (AMG) method for solving the (graph) Laplacian problem.  The corresponding linear system is defined in terms of the
following bilinear form:
\begin{eqnarray}\label{prob}
(Au,v) = 
\sum_{e\in\mathcal{E}} w_{e}\delta_e u\delta_e v + \sum_{i\in S_b}d_iu_iv_i
=(f,v),  
\end{eqnarray}
where $\mathcal{G}=(\mathcal{V},\mathcal{E})$ denotes an unweighted connected graph, $\mathcal{V}$
and $\mathcal{E}$ denote the set of vertices and edges of $\mathcal{G}$, respectively,  
and $\delta_e u = (u_i - u_j)$ for $e=(i,j)\in \mathcal{E}$. 
Note that the {\em lower-order terms}, $d_iu_i v_i, i\in S_b$, are included to account for various types of {\em boundary conditions}
for problems originating from discretization of partial differential equations (PDEs). 
If the lower-order terms are omitted and the wieghts $w_e = 1$, then the variational problem reduces to the graph Laplacian for a graph $\mathcal{G}$ that we focus on here.
The graph Laplacian, $A$, is then a symmetric and positive semi-definite matrix
and its kernel is the space spanned by the constant vector.

The main aim of the paper is to study the use of polynomial smoothing together with aggressive unsmoothed aggregation-based algebraic multigrid (UA-AMG) coarsening in developing an AMLI-cycle or k-cycle preconditioner 
~\cite{amli-1} for the graph Laplacian system.  
We consider the recently proposed polynomial based on the best approximation to $x^{-1}$ in the uniform norm~\cite{KVZ_09} in formulating the proposed UA-AMG algorithm.
A multilevel smoothed aggregation (SA) AMG algorithm using polynomial smoothers based on Chebychev approximations and its $V$-cycle convergence analysis are found in~\cite{PVanek_MBrezina_JMandel_2001a}.   We note that, these results are also used in~\cite{KVZ_09} to derive an SA two-level preconditioner with polynomial smoothing for diffusion problems.  In both methods,  the polynomial approximation is used to form (1) a smoother for the interpolation operator and (2) a relaxation scheme for the solver.  These preconditioners yield uniformly convergent methods provided polynomials of sufficiently large degree are used in both steps.  Further development and analysis of polynomial smoothers are found in~\cite{Adams03parallelmultigrid} and~\cite{BVV_2011,baker:2864}. 

Here, we consider an approach in which the polynomial smoother is used as the relaxation scheme in the AMG solver and interpolation is based on UA-AMG framework.  We show that using such plain aggregation based aggressive coarsening with a polynomial smoother in a AMLI cycle or k-cycle leads to a uniformly convergent method.  Generally, the use of unsmoothed (or plain) aggregation to construct the coarse space
without the use of interpolation smoothing has been observed to result in
slow convergence of a $V$-cycle multilevel iterative solver.
We note that recently it has been shown that plain aggregation-based coarsening approaches can lead to effective solvers for a variety of problems provided AMLI or k-cycles are used, e.g, such approaches have been developed and analyzed for the graph Laplacian in~\cite{LAMG_Report}, for more general $M$ matrices in~\cite{aggproved,brannick_local_stab_2011}, 
and for problems in quantum dynamics in~\cite{BrannickBrowerClarkOsborneRebbi_2007aa}. 
Generally, the use of AMLI cycles and UA-AMG typically leads to low grid and operator
complexities, limited fill-in in the coarse level operators, and reduces the arithmetic
complexity in the setup phase substantially.  The gains in the solve phase are often 
less pronounced since AMLI- and NAMLI-cycles use additional coarse-level corrections
to accelerate convergence of the UA-AMG method.

In Section \ref{matching}, we introduce a graph partitioning algorithm for constructing the coarse space.  Then, in Section \ref{matching2}, we establish an approximation property for such piecewise constant coarse spaces, which together with the stability estimates for such methods found in~\cite{brannick_local_stab_2011},  gives a spectral equivalence result that holds for the corresponding two-level method applied to graph Laplacian on general graphs.   The resulting estimate depends on the degree of the polynomial smoother and the coarsening ratio, i.e., the cardinality of the aggregates, and thus provides a way to adjust the polynomial degree to compensate for aggressive coarsening.  We note that the result is a special case of the general result found in~\cite{KVZ_09}. 
In the last section, we provide numerical experiments of the proposed multigrid approach applied to the graph Laplacian and show that the coarsening can be quite aggressive and still only a low degree polynomial is needed to obtain a scalable AMLI or k-cycle preconditioner.    

\section{Subspaces by graph partitioning and graph matching}
\label{matching}

We define a graph partitioning of $\mathcal{G}=(\mathcal{V},\mathcal{E})$ as a set of connected subgraphs 
$\mathcal{G}_{i}=(\mathcal{V}_{i},\mathcal{E}_{i})$ such that 
$\cup_{i}\mathcal{V}_{i}=\mathcal{V}
, \quad
\mathcal{V}_{i}\cap\mathcal{V}_{j} =\emptyset
, \quad
i\neq j.
$
In this paper, all subgraphs are assumed to be non empty and connected. 
The simplest non trivial example of such a graph partitioning is a matching, i.e, 
 a collection (subset $\mathcal{M}$) of edges in $\mathcal{E}$ such that no two
edges in $\mathcal{M}$ are incident. 
For a given graph partitioning, subspaces of $V=\mathbb{R}^{|\mathcal{V}|}$ are defined as
\begin{equation}\label{eq:coarse}
V_{H}=\{v \in V | \ v=\text{ constant on each $\mathcal{V}_{i}$ } \} .
\end{equation}
Note that each vertex in $\mathcal{G}$ corresponds to a connected subgraph $\mathcal{G}_i$ of $\mathcal{G}$
and every vertex of $\mathcal{G}$ belongs to exactly one such component. The vectors from
$V_H$ are constants on these connected subgraphs. 
The $\ell_2$ orthogonal projection on
$V_H$, which is denoted by $Q$, is defined as follows:
\begin{equation}\label{Q}
(Q v)_i = \frac{1}{|\mathcal{V}_{k}|}\sum_{j\in \mathcal{V}_{k}} v_j 
,\quad
\forall i \in \mathcal{V}_{k}. 
\end{equation}
Given a graph partitioning, 
the coarse graph $\mathcal{G}_{H}=\{\mathcal{V}_{H},\mathcal{E}_{H}\}$ is defined by 
assuming that all vertices in a subgraph form an equivalence class
and that $\mathcal{V}_{H}$ and $\mathcal{E}_{H}$ are the quotient set of $\mathcal{V}$ and $\mathcal{E}$ 
under this equivalence relation. 
That is, any vertex in $\mathcal{V}_{H}$ corresponds to a subgraph in the partitioning, 
and the edge $(i,j)$ exists in $\mathcal{E}_{H}$ if and only if 
the $i$-th and $j$-th subgraphs are connected in the graph $\mathcal{G}$. 

The algorithm we use in forming a graph partitioning is a variant of the approach we developed and tested for graphics processing units in~\cite{gpu}.  The 
procedure iteratively applies the following two steps:
\begin{enumerate}
\item[(A)]
Construct a set $S$ which contains coarse vertices by applying a maximal independent set algorithm to the graph of $A^k$.
\item[(B)]
Construct a subgraph for each vertex in $S$ by collecting vertices and edges of  the neighbors of vertices in $S$. 
\end{enumerate} 

\section{Two-level preconditioner with polynomial smoothing for the graph Laplacian}
\label{matching2}
A variational two-level method with one post smoothing step is defined as follows.  Given an approximation $w\in V$ to the solution $u$
of the graph Laplacian system, an update $v\in V$
is computed in two steps
\begin{enumerate}
\item $y = w + PA_H^{\dagger}P^T(f-Aw), \quad A_H = P^TAP$.
\item $v = y+R(f-Ay)$.
\end{enumerate}
We use $\dagger$ to denote the pseudo inverse of a matrix.
The corresponding error propagation operator of the two-level method is given 
by
\[
E_{TL} = (I-RA)(I-\pi_A), \quad \pi_A = PA_H^{\dagger} P^T A.
\]
Here, $E_{TL}$ is nonsymmetric and, thus, we consider
the following {\em symmetrization} to form the two-level preconditioner:
$
B=(I-E_{TL}E_{TL}^{*})A^{\dagger}, 
$
with
$^*$ denoting the adjoint with respect to the energy
inner product $(\cdot,\cdot)_A$.
We note that $|E_{TL}|^2_A = \rho(I-BA)$, where $\rho(X)$ is the spectral radius of the matrix $X$.  Further, 
if $\bar{R}$ satisfies  
$
\left( I - \bar{R} A \right) = \left( I - R A \right)^2$ so that 
$ \bar{R} = 2R - R A R$, then using that 
$\pi_A$ is an $A$-orthogonal projection 
on $\mbox{range}(P)$, it follows by direct computation that
$B = \bar{R} + \left( I - R A \right) P  A^{\dagger}_H P^T \left( I - A R \right).$

In~\cite{KVZ_09}, a spectral equivalence result for the preconditioner $B$ using 
a polynomial smoother based on the best approximation to $x^{-1}$ on a finite interval
$\left [\lambda_0,\; 
\lambda_1 \right ]$,
$0<\lambda_0 <\lambda_1$, in the
uniform norm ( $\|\cdot\|_{\infty}$) is derived.  
Here, $\lambda_0>0$ is any
lower bound on the spectrum of $A$
and $\lambda_1 = \| A \|_{\ell_\infty}$ is an approximation
to $\rho(A)$.  The unique solution to the minimization
problem 
\begin{equation}\label{1x-1-min}
q_m(x)= \arg\min\{ \|\frac1x-p\|_{\infty,[\lambda_0,\lambda_1]}, \quad p\in \mathcal{P}_{m}\}, 
\end{equation}
determines the polynomial approximation of degree $m$.
For details on the three-term 
recurrence used in its construction we refer to~\cite{KVZ_09}.  
Define 
\begin{eqnarray*}
E_m:= \max_{x\in [\lambda_0,\lambda_1]}|1-xq_m(x)|
=\max_{x\in [\lambda_0,\lambda_1]} x\cdot\left|\frac{1}{x}-q_m(x)\right|.
\end{eqnarray*}
Then, since $\lambda_1$ is a point of Chebyshev alternance
from
\cite[Theorem~2.1 and
Equation~(2.2)]{KVZ_09} for the error of approximation $E_m$ we have
\begin{eqnarray*}
E_m = \lambda_1
\left| \frac{1}{\lambda_1}- q_m(\lambda_1)\right|
= \left[\frac{2\lambda_1}{\lambda_1-\lambda_0}\right]\cdot \left[\frac{\delta^m}{a^2-1}\right] = \frac{2\kappa\delta^m}{(\kappa-1)(a^2-1)}.
\end{eqnarray*}
Here, we have denoted
$
\kappa=\frac{\lambda_1}{\lambda_0},$ 
$\delta = \frac{\sqrt{\kappa}-1}{\sqrt{\kappa}+1}$, and
$a = \frac{\kappa+1}{\kappa-1}.
$
Computing the error $E_m$ then gives
\[
E_m = \frac{\delta^m(\kappa-1)}{2}.
\]
A restriction on the degree $m$ is given by the requirement that
$q_m(\lambda_1)>0$. A sufficient condition for the
positivity of this polynomial (and
also necessary condition in many cases) is that
$\frac{1}{\lambda_1}-E_m >0$. Thus,  we need to find the smallest $m$
such that both $E_m< \rho$
and $q_m(\lambda_1)>0$. We then have that the polynomial is
positive if
\begin{eqnarray*}
\frac{\delta^m(\kappa-1)}{2} \le  \frac{1}{\lambda_1}
\quad\Rightarrow\quad
\delta^m\le \frac{2}{\lambda_1(\kappa-1)}.
\end{eqnarray*}
We note that from this it follows that 
$R=q_m(A)$ and hence $\bar{R}$ are
symmetric and positive definite, implying that the smoother 
in convergent in $A$-norm.
Also, to guarantee a damping
factor less than $\rho$ on the interval $[\lambda_0,\lambda_1]$, we have
\begin{eqnarray*}
\frac{\delta^m(\kappa-1)}{2} \le  \rho
\quad\Rightarrow\quad
\delta^m\le \frac{2\rho}{\kappa-1}.
\end{eqnarray*}
Thus, the minimal $m$ that guarantees both properties are satisfied is given by
\begin{equation}\label{eq:bound-on-m}
m \ge
\frac{1}{|\log\delta|}\max\left\{ \left|\log \frac{2\rho}{\kappa-1}\right|,
\left|\log \frac{2}{\lambda_1(\kappa-1)}\right|
\right\}.
\end{equation}   
The spectral equivalence result that we adopt to analyze a two-level method based on plain aggregation  with this polynomial smoother follows from this smoothing estimate and the assumptions 
of stability and an approximation property of the coarse space 
$V_H$:
for any $v\in V_h$, 
\begin{equation}\label{eqn:approx}
c_p^{-1}\|v-Qv\|^2 +|v-Qv|^2_{A} \le c_1 |v|^2_A,
\end{equation}
where $|\cdot|_A$ denotes the $A$ semi norm.
Recall that in this paper $Q$ is the $\ell_2$ projection on the span of
$\{1_l\}_{l=1}^{n_H}$ (p.w. constant projection) and, thus, 
this inequality holds also in the case that $v$ is in the kernel of $A$, since 
then, all the terms are equal to zero.  

Assume that $V_H$ is such that the above approximation and stability assumptions hold and the polynomial $q_m$ is chosen such that (\ref{eq:bound-on-m}) holds for a fixed value $\lambda_0$. Then, the following spectral equivalence holds 
\begin{equation}
v^TAv \leq v^TB^{\dagger}v \leq K_{TG}\;
v^TAv,\quad K_{TG}=8+ 8c_{1}\left[ c_{n_z} c
_p c_s + 1\right]. 
\end{equation}
This result is a special case of Theorem 4.6 in~\cite{KVZ_09}, refined for 
unsmoothed aggregation applied to the graph Laplacian.
Here, $c_{n_z}$ is a constant that depends on the number of nonzeros per row of $A$, the constant $c_1$ involves the stability of $Q$ in $A$-norm and the constant $c_p$ arises from the weak approximation property, and as we show below, 
depends on the cardinality and the diameter of the subgraphs in the graph partitioning.  The constant $c_s = \frac{ \ln^2 m}{m^2}$, where $m$ is the degree of the polynomial.  
Thus, given a partitioning of the fine-level graph into subgraphs, $\mathcal{G} = \cup_{l=1}^{n_H} \mathcal{G}_l$, it is possible to choose the degree of the polynomial $m$ sufficiently large to control
the constant $c_p$ and hence $K_{TG}$ in the above spectral equivalence estimate.  
This result is derived from the following estimate (see Corollary 4.4 in~\cite{KVZ_09})
\begin{equation}\label{eq:xz}
v^TB^{\dagger} v \leq 4 \; \inf_{v_h \in V_H} \bigg[|v_H|_A^2 +  \lambda c_s
\|v - v_H\|^2 + |v - v_H|_A^2\bigg].
\end{equation}
A similar result for smoothed aggregation based on Chebyshev polynomial approximations is found in~\cite{BVV_2011}.  

Next, we establish the approximation property for the p.w. constant coarse space $V_H$ as defined in (\ref{eq:coarse}) for the graph Laplacian.  Suppose that $\mathcal{V} = \{1,\ldots,n\}$ is partitioned into nonoverlapping
subsets: $\mathcal{V}=\cup_{l=1}^{n_H} \mathcal{V}_l, n_l=|\mathcal{V}_l|.$
Each set of vertices defines a subgraph $\mathcal{G}_\ell$ whose vertex
set is $\mathcal{V}_l$ and whose edges $\mathcal{E}_l$ are a subset of
$\mathcal{E}$, where $(i,j)\in \mathcal{E}_l$ if and only if 
\emph{both} $i$ and $j$ are in $\mathcal{V}_l$. Denote the graph Laplacian
associated with the subgraph $\mathcal{G}_l$ by $A_l$. 
Let $1$ denote the constant vector on $\mathcal{V}$ and $1_l$ the constant 
vector on $\mathcal{V}_l$ extended by $0$ outside $\mathcal{V}_l$. Let $\lambda_l$ 
be the smallest positive eigenvalue of the graph Laplacian on $\mathcal{G}_l$, namely, 
$\lambda_l$ is defined as $\displaystyle{
\lambda_l = \min_{v:\;(v,1_l)=0} \frac{(A_l v,v)}{\|v\|^2}}.
$
Here, the minimum is taken over all  $v\in \mathbb{R}^{n_l}$. Given $v\in \mathbb{R}^n$ define $\|v\|_{\mathcal{G}_l}^2 = \sum_{j\in \mathcal{V}^l} v_j^2$, which is the $\ell_2$ norm on the subgraph $\mathcal{G}_l$. 
Now, since
$((v-Qv),1_l)=0$,
we have  
$
\|v-Qv\|^2_{\mathcal{G}_l} \le \lambda_l^{-1}\sum_{e\in
  \mathcal{E}_l}(\delta_e v)^2.
$ 
Thus, 
\begin{eqnarray}\label{eq:wap}
\|v-Qv\|^2 = 
\sum_{l=1}^{n_c}\|v-Q v\|^2_{\mathcal{G}_l} \le 
\sum_{l=1}^{n_c}\lambda_l^{-1 }\sum_{e\in\mathcal{E}_l} (\delta_e v)^2
 \le  c_p\sum_{e\in\mathcal{E}} (\delta_e v)^2 = c_p (Av,v). 
\end{eqnarray}
The last step follows from the definition of $c_p$
and the observation that since $\cup_{l}
\mathcal{E_\ell}\subset\mathcal{E}$,we have that for any $v\in \mathbb{R}^n$, 
$
\sum_{l=1}^{n_H}\sum_{e\in \mathcal{E}_l} (\delta_e v)^2 \le 
\sum_{e\in \mathcal{E}} (\delta_e v)^2 =(Av,v). 
$
Note that this latter result holds since the second sum is over a larger set. 
For {\em shape regular} subgraphs, $\mathcal{G}_l$, the local constants $\lambda_l^{-1}$ can be bounded in terms of $|\mathcal{V}_l| \cdot \mbox{diam}(\mathcal{G}_l)$ using Cheeger's inequality~\cite{Cheeger}.  Here, $\mbox{diam}(\mathcal{G}_l)$ denotes the diameter of the longest path in the $l$th subgraph.  
A similar technique is considered in~\cite{aggproved}, in which the constants $\lambda_\ell^{-1}$ are computed by solving local eigenvalue problems.

In~\cite{BCKZ2012}, commuting relations involving a certain projection, $\Pi$, the p.w. constant projection $Q$, the discrete gradient operator, $B$, and $B^T$ on the graph, $\mathcal{G}$, are introduced and are then used to derive a stability estimate of the form
$$
|Q|^2_{A}
=
\sup_{v:\;(v,1)=0}\frac{|Qv|^2_{A}}{|v|^2_{A}} \leq \|\Pi\|^2 \leq c_0\; ,
$$
where $c_0$ is a constant that
depends on the shape and alignment of the subgraphs, 
but not on the dimension, $|\mathcal{V}|$, of the graph Laplacian, $A$.  
It is noteworthy that this bound holds for general graphs with few assumptions and, further, 
that, since $\Pi$ is constructed one row at a time, this estimate allows 
local energy estimates that can be used in forming the graph
partitioning.  A similar approach was considered in
\cite{LAMG_Report}.

Given the above approximation and stability estimates 
and using that $|v_H|_A \leq c_0 |v|_A$, $v_H = Qv$, it follows that 
the inequality in (\ref{eqn:approx}) holds with 
$c_1 = 2c_0+3$ and $c_p$ given in 
(\ref{eq:wap}).  This, in turn, implies the spectral equivalence of the 
two-level preconditioner based on a p.w. constant coarse space
$V_H$ for the graph Laplacian.
We remark that the Galerkin coarse-level operator $A_H = 
P^TAP$ is generally a weighted graph Laplacian of the form
$A_H = B_H^T D B_H$, where $D$ is a diagonal weight matrix with strictly positive entries and $B_H$ is the discrete gradient operator defined on the coarse graph $\mathcal{G}_H(\mathcal{V}_H,\mathcal{E}_H)$.  Similar stability and approximation
properties of piece-wise constant coarse spaces can be established 
in this more general setting as well and, then, a similar proof
of the spectral equivalence result follows with minor modifications.
Alternatively, it is possible to replace the weighted graph Laplacian
with an unweighted one on the same graph and derive a spectral
equivalence result between the two.  The latter result, in turn, 
again can be used to establish a spectral equivalence result for this
modified two-level method.

\section{Numerical results}\label{numerics}

We apply the proposed aggregation based
preconditioner to graph Laplacians resulting from finite element discretizations of the scalar Laplace problem.  
We consider both stationary AMLI-cycle and N-AMLI-cycle (k-cycle) preconditioners.  
For details on the theory and the implementation of the AMLI and N-AMLI methods 
we refer to~\cite{amli-1}.  In the AMLI approach, we use the polynomial based on the best approximation to $x^{-1}$ in the uniform norm to form a the preconditioner between any two successive levels of the multilevel hierarchy, see~\cite{KVZ_09}.  In the N-AMLI-cycle, a nonlinear PCG (NPCG) method is applied recursively to solve the coarse-level equations.   
The AMLI-cycle is used as a preconditioner for the CG method on the finest level and the N-AMLI-cycle is applied as a preconditioner to the NPCG iteration.  To limit the memory requirements of the nonlinear scheme we restart the outer fine-level NPCG method 
every five iterations.

In all tests, the maximal independent set algorithm used in the aggregation process for constructing the coarse spaces is applied to the graph of $A^4$, yielding a coarsening factor of roughly $n/n_H = 30$ between any two successive levels.  The problem is coarsened until the size of the coarsest level is less than 100.  As the relaxation method in the multilevel solver we use the polynomial smoother based on the best approximation to $x^{-1}$ on the interval $\left [\lambda_0,\; 
\lambda_1 \right ]$, where the estimate of the largest eigenvalue is computed as
$\lambda_1 = \| A \|_{\ell_\infty}$ and we set $\lambda_0=\lambda_1/10$.  Thus, taking the degree as $m=4$
in the polynomial smoother ensures the inequality (\ref{eq:bound-on-m}) holds.  
We test $W$-cycle AMLI and N-AMLI preconditioners with such smoother.  
The stopping criteria for the flexible preconditioned conjugate gradient iteration is set to a $10^{-8}$ reduction in 
the relative $A$ norm of the error and the number of iterations needed to reach this
tolerance in the different tests are reported.   

In Table \ref{tab:cf30}, we report results of the proposed method for graph Laplacians arising from discretizing the Poisson problem on structured and unstructured
meshes.  We compare the performance of a stationary AMLI with a N-AMLI, both using the same multilevel hierarchy obtained by applying the aggregation algorithm to the same Poisson problem with Neumann boundary conditions discretized using standard linear Finite Elements.   For the structured meshes we consider a 2d unit square domain with $n^2$ unknowns (left) and a 3d unit cube 
domain with $n^3$ unknowns (middle).  
Results for more general graphs (right), coming from unstructured meshes resulting from triangulations of the 3d unit cube, are also included. 
The unstructured mesh is formed by adding a random vector 
of length $h/2$, 
where $h$ is the grid length, 
to each vertex of a structured triangulation, 
followed by a Delaunay triangulation. 
The (AMLI) N-AMLI method yields a (nearly) scalable solver with low grid and operator complexities -- in all tests the grid complexities $\frac{\sum_{j=0}^{J} n_ j}{n_0}$ were less than
1.03 and the operator complexities $\frac{\sum_{j=0}^{J} nnz(A_ j)}{nnnz(A_0)}$ were less than 1.04.  
\begin{table}\label{tab:cf30}
\medskip
\centerline{
\begin{tabular}{ccc}
2d struct. & 3d struct. & 3d unstruct. \\
\begin{tabular}{|c|c|c|c|c|}
\hline
$n$ & AMLI & N-AMLI \\
\hline
 $512^{2}$  & 20  & 19 \\
\hline
 $1024^{2}$  & 22  & 20 \\
\hline
 $2048^{2}$  & 23  & 21 \\
\hline
 $4096^{2}$  & 24  & 21 \\
\hline
\end{tabular}
\hspace{.4cm}
&
\begin{tabular}{|c|c|c|c|c|}
\hline
$n$ & AMLI & N-AMLI \\
\hline
 $32^{3}$  & 22  & 20 \\
\hline
 $64^{3}$  & 23  & 22 \\
\hline
 $128^{3}$&  23 & 22  \\
\hline
 $256^{3}$  & 25  & 22 \\
\hline
\end{tabular}
\hspace{.4cm}
&
\begin{tabular}{|c|c|c|c|c|}
\hline
$n$ & AMLI & N-AMLI \\
\hline
 $32^{3}$  & 24  & 21 \\
\hline
 $64^{3}$  & 25  & 23 \\
\hline
 $128^{3}$&  27 & 24  \\
\hline
 $256^{3}$  & 28  & 24 \\
\hline
\end{tabular}
\end{tabular}}
\caption{ Results of $W(1,1)$ AMLI and nonlinear AMLI preconditioners with degree $m=4$ polynomial smoother for the Poisson problem.}
\end{table}

\section{Conclusion}\label{conclusion}
An algebraic graph partitioning algorithm for aggressive coarsening is developed and a two-level convergence theory of the resulting solver with polynomial smoother is developed.  It is shown numerically that the resulting N-AMLI approach with polynomial smoother yields an efficient solver for graph Laplacian problems coming from
Finite Element discretizations of the Poisson problem.  The graph partitioning algorithm, intended  for unweighted graphs, is designed to select shape regular aggregates of arbitrary size and, thus, can be used to obtain predefined coarsening factors.  The use of an unsmoothed
aggregation form of aggressive coarsening results in low overall grid and operator complexities
and limited fill-in in the coarse-level operators.  It further significantly simplifies the triple 
matrix product to simple summations of entries of the fine-level matrix.  
\bibliographystyle{plain}
\bibliography{brannickj_mini_13_biblio}

\providecommand{\noopsort}[1]{}
\begin{thebibliography}{10}

\bibitem{Adams03parallelmultigrid}
Mark Adams, Marian Brezina, Jonathan Hu, and Ray Tuminaro.
\newblock Parallel multigrid smoothing: polynomial versus gauss-seidel.
\newblock {\em J. Comp. Phys}, 188:593--610, 2003.

\bibitem{amli-1}
O.~Axelsson and P.~S. Vassilevski.
\newblock Algebraic multilevel preconditioning methods. {I}.
\newblock {\em Numer. Math.}, 56(2-3):157--177, 1989.

\bibitem{baker:2864}
Allison~H. Baker, Robert~D. Falgout, Tzanio~V. Kolev, and Ulrike~Meier Yang.
\newblock Multigrid smoothers for ultraparallel computing.
\newblock {\em SIAM Journal on Scientific Computing}, 33(5):2864--2887, 2011.

\bibitem{BrannickBrowerClarkOsborneRebbi_2007aa}
J.~Brannick, R.~Brower, M~Clark, J~Osborne, and C.~Rebbi.
\newblock Adaptive multigrid for lattice {QCD}.
\newblock {\em Phys. Rev. Let.}, 100(4):041601, 2008.

\bibitem{gpu}
J.~Brannick, Y.~Chen, X.~Hu, and L.~Zikatanov.
\newblock Parallel unsmoothed aggregation algebraic multigrid algorithms on
  gpus.
\newblock In {\em Numerical Solution of Partial Differential Equations: Theory,
  Algorithms and their Applications, Springer Series in Mathematics and
  Statistics}. Springer, Berlin/Heidelberg, Accepted 2012.

\bibitem{BCKZ2012}
J.~Brannick, Y.~Chen, J.~Krauss, and L.~Zikatanov.
\newblock An algebraic multigrid method based on matching of graphs.
\newblock {\em Lecture notes in Computational Science and Engineering}, 2012.

\bibitem{brannick_local_stab_2011}
J.~Brannick, Y.~Chen, and L.~Zikatanov.
\newblock An algebraic multilevel method for anisotropic elliptic equations
  based on subgraph matching.
\newblock {\em Numer. Linear Algebra Appl.}, 19:279--295, 2012.

\bibitem{BVV_2011}
M.~Brezina, P.Van\v{e}k, and P.~Vassilevski.
\newblock An improved convergence analysis of smoothed aggregation algebraic
  multigrid.
\newblock {\em Journal of Numerical Linear Algebra and Applications}, 19:pages
  441--469, 2012.

\bibitem{Cheeger}
J.~Cheeeger.
\newblock A lower bound for the smallest eigenvalue of the {L}aplacian.
\newblock In Robert~C. Gunning, editor, {\em Problems in Analysis, A Symposium
  in Honor of Salomon Bochner}, pages 195--199. Princeton Univ. Press, New
  Jersey, 1970.

\bibitem{KVZ_09}
J.~Kraus, P.~Vassilevski, and L.~Zikatanov.
\newblock Polynomial of best uniform approximation to $x^{-1}$ and smoothing in
  two-level methods.
\newblock {\em Comput. Methods Appl. Math.}, 12:448--468, 2012.
\newblock Also available at arXiv:1002.1859v3 [math.NA].

\bibitem{LAMG_Report}
O.E. Livne and A.~Brandt.
\newblock Lean algebraic multigrid ({LAMG}): Fast graph {L}aplacian linear
  solver.
\newblock {\em SIAM Journal of Scientific Computing}, 34:B499--B522, 2012.

\bibitem{aggproved}
A.~Napov and Y.~Notay.
\newblock An algebraic multigrid method with guaranteed convergence rate.
\newblock {\em SIAM J. Sci. Comput.}, 34:A1079ÐA1109, 2012.

\bibitem{PVanek_MBrezina_JMandel_2001a}
P.~Van{\v e}k, M.~Brezina, and J.~Mandel.
\newblock Convergence of algebraic multigrid based on smoothed aggregation.
\newblock {\em Numer. Math.}, 88:559--579, 2001.

\end{thebibliography}
\end{document}